\documentclass[12pt]{article}
\usepackage{amssymb,graphicx}
\newcommand{\wt}{\widetilde}
\newcommand{\R}{\mathbb R}
\newcommand{\C}{\mathbb C}

\begin{document}
\thispagestyle{empty}
\hfill
\vskip 3.8truecm
\centerline{{\large The Beckman-Quarles theorem for continuous}}
\centerline{{\large mappings from ${\R}^2$ to ${\C}^2$}}
\vskip 0.8truecm
\centerline{{\large Apoloniusz Tyszka}}
\vskip 0.8truecm
\rightline{IMUJ Preprint 2002/   \hspace{2.4cm}}
\vskip 1.8truecm
\begin{abstract}
Let $\varphi:{\C}^2 \times {\C}^2 \to \C$,
$\varphi((x_1,x_2),(y_1,y_2))=(x_1-y_1)^2+(x_2-y_2)^2$.
We say that $f:{\R}^2 \to {\C}^2$ preserves distance $d~\geq~0$ if for
each $x,y \in {\R}^2$ $\varphi(x,y)=d^2$ implies $\varphi(f(x),f(y))=d^2$.
We prove that if $x,y \in {\R}^2$ and
$|x-y|=(2\sqrt{2}/3)^k \cdot (\sqrt{3})^l$ ($k,l$ are non-negative
integers) then there exists a finite set
$\{x,y\} \subseteq S_{xy} \subseteq {\R}^2$ such that each
unit-distance preserving mapping from
$S_{xy}$ to ${\C}^2$ preserves the distance
between $x$ and $y$. It implies that each continuous map from
${\R}^2$ to ${\C}^2$ preserving unit distance
preserves all distances.
\end{abstract}
\vskip 0.2truecm
\normalsize
\par
The classical Beckman-Quarles theorem states that each unit-distance
preserving mapping from ${\R}^n$ to ${\R}^n$ ($n \geq 2$) is an isometry,
see \cite{Beckman}, \cite{Benz}, \cite{Everling} and~\cite{Lester}.
Author's discrete form of this theorem
(\cite{Aequationes},\cite{Yerevan}) states that if
$x,y \in {\R}^n$ ($n \geq 2$) and $|x-y|$ is an algebraic number then
there exists a finite set $\{x,y\} \subseteq S_{xy} \subseteq {\R}^n$
such that each unit-distance preserving mapping from $S_{xy}$ to
${\R}^n$ preserves the distance between $x$ and $y$.
\vskip 0.2truecm
\par
Let $\varphi_n:{\C}^n \times {\C}^n \to \C$,
$\varphi_n((x_1,...,x_n),(y_1,...,y_n))=(x_1-y_1)^2+...+(x_n-y_n)^2$.
We say that $f:{\R}^n \to {\C}^n$ preserves distance $d \geq 0$ if for
each $x,y \in {\R}^n$ $\varphi_n(x,y)=d^2$ implies $\varphi_n(f(x),f(y))=d^2$.
In \cite{zlozona} the author proved that each continuous mapping
from ${\R}^n$ to ${\C}^n$ ($n \geq 3$) preserving unit distance preserves
all distances. In this paper we prove it for $n=2$, similarly to the
case $n \geq 3$ the proof is based on calculations using the
Cayley-Menger
determinant.
\vskip 0.2truecm
\par
{\bf Proposition~1} (\cite{zlozona}, cf. \cite{Blumenthal}, \cite{Borsuk}).
The points $c_{1}=(z_{1,1},...,z_{1,n}),...,
c_{n+1}=(z_{n+1,1},...,z_{n+1,n}) \in {\C}^n$ are affinely dependent
if and only if their Cayley-Menger determinant
$$
\det \left[
\begin{array}{ccccc}
 0  &  1                       &  1                       & ... & 1                         \\
 1  & \varphi_n(c_{1},c_{1})   & \varphi_n(c_{1},c_{2})   & ... & \varphi_n(c_{1},c_{n+1})  \\
 1  & \varphi_n(c_{2},c_{1})   & \varphi_n(c_{2},c_{2})   & ... & \varphi_n(c_{2},c_{n+1})  \\
... & ...                      & ...  	                  & ... & ...                       \\
 1  & \varphi_n(c_{n+1},c_{1}) & \varphi_n(c_{n+1},c_{2}) & ... & \varphi_n(c_{n+1},c_{n+1})\\
\end{array}\;\right]
$$
\par
\noindent
equals $0$.
\par
\vskip 0.2truecm
{\it Proof.} It follows from the equality
$$
\left(
\det \left[
\begin{array}{ccccc}
z_{1,1}   & z_{1,2}   & ... &  z_{1,n}  & 1  \\
z_{2,1}   & z_{2,2}   & ... &  z_{2,n}  & 1  \\
  ...     &  ...      & ... &  ...      & ...\\
z_{n+1,1} & z_{n+1,2} & ... & z_{n+1,n} & 1  \\
\end{array}
\right] \right)^2=
$$
$$
\frac{(-1)^{n+1}}{2^{n}} \cdot
\det \left[
\begin{array}{ccccc}
 0  &  1                       & 1                        & ... &  1                        \\
 1  & \varphi_n(c_{1},c_{1})   & \varphi_n(c_{1},c_{2})   & ... & \varphi_n(c_{1},c_{n+1})  \\
 1  & \varphi_n(c_{2},c_{1})   & \varphi_n(c_{2},c_{2})   & ... & \varphi_n(c_{2},c_{n+1})  \\
... & ...                      & ...                      & ... & ...                       \\
 1  & \varphi_n(c_{n+1},c_{1}) & \varphi_n(c_{n+1},c_{2}) & ... & \varphi_n(c_{n+1},c_{n+1})\\
\end{array}
\right].
$$

\vskip 0.2truecm
\par
{\bf Proposition~2} (see \cite{zlozona} for $k=2$, cf. \cite{Blumenthal}, \cite{Borsuk}).
For each points $c_{1},...,c_{n+k} \in {\C}^n$ ($k=2,3,4,...$) their
Cayley-Menger determinant equals $0$ i.e.
$$
\det \left[
\begin{array}{ccccc}
 0  &  1                       &  1                       & ... & 1                         \\
 1  & \varphi_n(c_{1},c_{1})   & \varphi_n(c_{1},c_{2})   & ... & \varphi_n(c_{1},c_{n+k})  \\
 1  & \varphi_n(c_{2},c_{1})   & \varphi_n(c_{2},c_{2})   & ... & \varphi_n(c_{2},c_{n+k})  \\
... & ...                      & ...  	                  & ... & ...                       \\
 1  & \varphi_n(c_{n+k},c_{1}) & \varphi_n(c_{n+k},c_{2}) & ... & \varphi_n(c_{n+k},c_{n+k})\\
\end{array}\;\right]
={\rm 0}.
$$
\vskip 0.3truecm
{\it Proof.} Assume that
$c_{1}=(z_{1,1},...,z_{1,n}),...,c_{n+k}=(z_{n+k,1},...,z_{n+k,n})$.
The points $\wt{c}_{1}=(z_{1,1},...,z_{1,n},0,...,0)$,
$\wt{c}_{2}=(z_{2,1},...,z_{2,n},0,...,0)$, . . . ,
$\wt{c}_{n+k}=(z_{n+k,1},...,z_{n+k,n},0,...,0) \in
{\C}^{n+k-1}$ are affinely dependent.
Since $\varphi_n(c_i,c_j)=\varphi_{n+k-1}(\wt{c}_{i},\wt{c}_{j})$
$(1 \leq i \leq j \leq n+k)$ the Cayley-Menger determinant
of points $c_{1},...,c_{n+k}$ is equal to the Cayley-Menger determinant
of points $\wt{c}_{1},...,\wt{c}_{n+k}$ which equals $0$ according to
Proposition~1.
\vskip 0.2truecm
\par
From Proposition~1 we obtain the following Propositions~3a and~3b.
\vskip 0.2truecm
\par
{\bf Proposition~3a.} If $c_{1},c_{2},c_{3} \in {\C}^2$
and $\varphi_2(c_1,c_2)=\varphi_2(c_1,c_3)=\varphi_2(c_2,c_3) \in (0,\infty)$ ,
then $c_{1}$, $c_{2}$, $c_{3}$ are affinely independent.
\vskip 0.2truecm
\par
{\bf Proposition~3b.} If $d>0$, $c_{1}$, $c_{2}$, $c_{3}$ $\in$ ${\C}^2$ and
$\varphi_2(c_1,c_2)=2d^2$, $\varphi_2(c_1,c_3)=3d^2$, 
$\varphi_2(c_2,c_3)=9d^2$, then $c_{1}$, $c_{2}$, $c_{3}$ are
affinely independent.
\vskip 0.3truecm
\par
{\bf Proposition~4} (see \cite{zlozona} for ${\C}^n$, cf. \cite{Borsuk} p. 127 for ${\R}^n$).
If $x$, $y$, $c_{0}$, $c_{1}$, $c_{2}$ $\in$ ${\C}^2$, 
$\varphi_2(x,c_{0})=\varphi_2(y,c_{0})$, 
$\varphi_2(x,c_{1})=\varphi_2(y,c_{1})$,
$\varphi_2(x,c_{2})=\varphi_2(y,c_{2})$
and $c_0$, $c_1$, $c_2$ are affinely independent,
then $x=y$.
\vskip 0.2truecm
\par
{\it Proof.} Computing we obtain that the vector
$\overrightarrow{xy}:=[s_1,s_2]$ is perpendicular
to each of the linearly independent vectors
$\overrightarrow{c_{0}c_{1}}$, $\overrightarrow{c_{0}c_{2}}$.
Thus the vector $\overrightarrow{xy}$ is perpendicular
to every linear combination of vectors 
$\overrightarrow{c_{0}c_{1}}$ and $\overrightarrow{c_{0}c_{2}}$.
In particular, the vector
$\overrightarrow{xy}=[s_1,s_2]$ is perpendicular
to the vector $[\bar{s_1},\bar{s_2}]$, where $\bar{s_1},\bar{s_2}$
denote numbers conjugate to the numbers $s_1,s_2$. Therefore
$\overrightarrow{xy}=0$ and the proof is complete.
\vskip 0.2truecm
\par
From Kronecker's theorem (see \cite{Hardy}) we obtain
the next proposition.
\vskip 0.2truecm
\par
{\bf Proposition~5.} The set
$\{(2\sqrt{2}/3)^k \cdot (\sqrt{3})^l: k,l \in \{0,1,2,...\}\}$ is a
dense subset of $(0,\infty)$.
\vskip 0.2truecm
\par
\vskip 0.3truecm
\par
{\bf Theorem~1.} If $x,y \in {\R}^2$ and
$|x-y|=(2\sqrt{2}/3)^k \cdot (\sqrt{3})^l$ ($k,l$ are non-negative
integers), then
there exists a finite set $\{x,y\} \subseteq S_{xy} \subseteq {\R}^2$
such that each unit-distance preserving mapping from $S_{xy}$ to ${\C}^2$
preserves the distance between $x$ and $y$.
\vskip 0.3truecm
\par
{\it Proof.}
Let $D$ denote the set of all positive numbers $d$
with the following property:
\begin{description}
\item{$(\ast)$}
if $x,y \in {\R}^2$ and $|x-y|=d$ then there exists a finite set
$\{x,y\} \subseteq S_{xy} \subseteq {\R}^2$ such that any map
$f:S_{xy}\rightarrow {\C}^2$ that preserves unit distance
preserves also the distance between $x$ and $y$.
\end{description}
\vskip 0.2truecm
\par
{\bf Lemma~1} (see \cite{zlozona} for $n \geq 3$).
If $d \in D$ then $\sqrt{3} \cdot d \in D$.
\vskip 0.2truecm
\par
{\it Proof.}
Let $ d \in D$, $x,y \in {\R}^2$ and $|x-y|=\sqrt{3} \cdot d$.
Using the notation of Figure~1 we show that
$$
S_{xy}:=S_{y\wt{y}} \cup \bigcup_{i=1}^2 S_{xp_i} \cup \bigcup_{i=1}^2 S_{yp_i} \cup S_{p_1p_2}
\cup
\bigcup_{i=1}^2 S_{x\wt{p_i}}
\cup
\bigcup_{i=1}^2 S_{\wt{y}\wt{p_i}}
\cup S_{\wt{p_1}\wt{p_2}}
$$
satisfies condition $(\ast)$.
\\
\includegraphics[width=137mm,height=70mm]{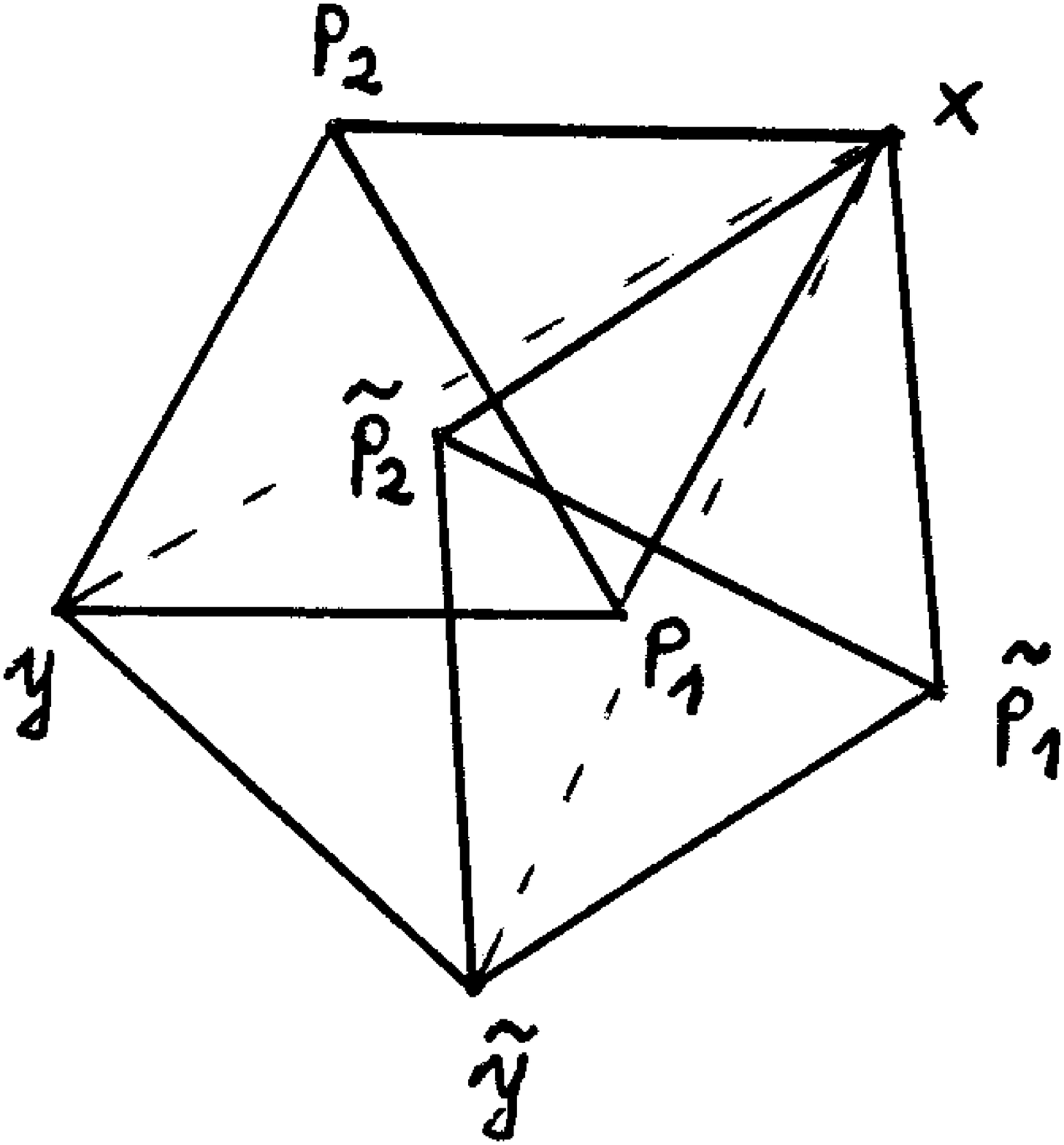}
\\
\centerline{Figure~1}
\centerline{$|x-y|=|x-\wt{y}|=\sqrt{3} \cdot d$, $|y-\wt{y}|=|p_1-p_2|=|\wt{p_1}-\wt{p_2}|=d$}
\centerline{$|x-p_i|=|y-p_i|=|x-\wt{p_i}|=|\wt{y}-\wt{p_i}|=d$ ($i=1,2$)} 
\\
\normalsize
\\
Assume that $f: S_{xy} \rightarrow {\C}^2$
preserves unit distance. Since
$$
S_{xy} \supseteq
S_{y\wt{y}}
\cup
\bigcup_{i=1}^{2}S_{xp_{i}}
\cup
\bigcup_ {i=1}^{2}S_{yp_{i}}
\cup
S_{p_{1}p_{2}}
$$
\noindent
we conclude that $f$ preserves the distances between
$y$ and $\wt{y}$,
$x$ and $p_{i}$ ($i=1,2$), $y$ and $p_{i}$ ($i=1,2$),
$p_{1}$ and $p_{2}$.
Hence $\varphi_2(f(y),f(\wt{y}))=
\varphi_2(f(x),f(p_i))=\varphi_2(f(y),f(p_i))=\varphi_2(f(p_1),f(p_2))=d^2$
($i=1,2$).
By Proposition~2 the Cayley-Menger determinant of points
$f(x)$, $f(p_{1})$, $f(p_{2})$, $f(y)$ equals $0$ i.e.
\footnotesize
$$
\det \left[
\begin{array}{ccccccc}
0 &  1                     &  1                       &  1                       &  1                    \\
1 & \varphi_2(f(x),f(x))   & \varphi_2(f(x),f(p_1))   & \varphi_2(f(x),f(p_2))   & \varphi_2(f(x),f(y))  \\
1 & \varphi_2(f(p_1),f(x)) & \varphi_2(f(p_1),f(p_1)) & \varphi_2(f(p_1),f(p_2)) & \varphi_2(f(p_1),f(y))\\
1 & \varphi_2(f(p_2),f(x)) & \varphi_2(f(p_2),f(p_1)) & \varphi_2(f(p_2),f(p_2)) & \varphi_2(f(p_2),f(y))\\
1 & \varphi_2(f(y),f(x))   & \varphi_2(f(y),f(p_1))   & \varphi_2(f(y),f(p_2))   & \varphi_2(f(y),f(y))  \\
\end{array}
\right]
=0.
$$
\normalsize
Denoting $t=\varphi_2(f(x),f(y))$ we obtain
$$
\det \left[
\begin{array}{ccccccc}
0 &  1  &  1  &  1  &  1 \\
1 &  0  & d^2 & d^2 &  t \\
1 & d^2 &  0  & d^2 & d^2\\
1 & d^2 & d^2 &  0  & d^2\\
1 &  t  & d^2 & d^2 &  0 \\
\end{array}
\right]
=0.
$$
\par
\noindent
Computing this determinant we obtain
$$
2d^{2}t \cdot (3d^2-t)=0.$$
Therefore
$$t=\varphi_2(f(x),f(y))=\varphi_2(f(y),f(x))=(\sqrt{3} \cdot d)^2$$
or $$t=\varphi_2(f(x),f(y))=\varphi_2(f(y),f(x))=0.$$
Analogously we may prove that
$$\varphi_2(f(x),f(\wt{y}))=\varphi_2(f(\wt{y}),f(x))=
(\sqrt{3} \cdot d)^2$$
or $$\varphi_2(f(x),f(\wt{y}))=\varphi_2(f(\wt{y}),f(x))=0.$$
If $t=0$ then the points $f(x)$ and $f(y)$ satisfy:
\\
\\
\centerline{$\varphi_2(f(x),f(x))=0=\varphi_2(f(y),f(x))$,}
\\
\centerline{$\varphi_2(f(x),f(p_1))=d^2=\varphi_2(f(y),f(p_1))$,}
\\
\centerline{$\varphi_2(f(x),f(p_2))=d^2=\varphi_2(f(y),f(p_2))$.}
\\
\par
\noindent
By Proposition~3a the points $f(x)$, $f(p_1)$, $f(p_2)$ are affinely
independent. Therefore by Proposition~4 $f(x)=f(y)$ and consequently
$$
d^2=\varphi_2(f(y),f(\wt{y}))=\varphi_2(f(x),f(\wt{y}))
\in \{(\sqrt{3} \cdot d)^2,~0 \}.
$$
Since $d^2 \neq (\sqrt{3} \cdot d)^2$ and $d^2 \neq 0$ we conclude that
the case $t=0$ cannot occur. This completes the proof of Lemma~1.
\vskip 0.5truecm
\par
{\bf Lemma~2.} If $d \in D$ then $3 \cdot d \in D$.
\vskip 0.2truecm
\par
{\it Proof.} It follows from Lemma~1 because
$3 \cdot d=\sqrt{3} \cdot (\sqrt{3} \cdot d)$.
\vskip 0.2truecm
\par
{\bf Lemma~3.} If $d \in D$ then $2 \cdot d \in D$.
\vskip 0.2truecm
\par
{\it Proof.}
Let $d \in D$, $x,y \in {\R}^2$ and $|x-y|=2 \cdot d$.
Using the notation of Figure~2 we show that
\vskip 0.1truecm
\centerline{$S_{xy}:=
\bigcup \{S_{ab}:a,b\in \{x,y,p_{1},p_{2},p_{3}\},
|a-b|=d \vee |a-b|=\sqrt{3} \cdot d\}$}
\vskip 0.1truecm
\par
\noindent
(where $S_{xp_3}$ and $S_{yp_2}$ are known to exist by Lemma~1)
satisfies condition $(\ast)$.
\\
\includegraphics[width=137mm,height=63mm]{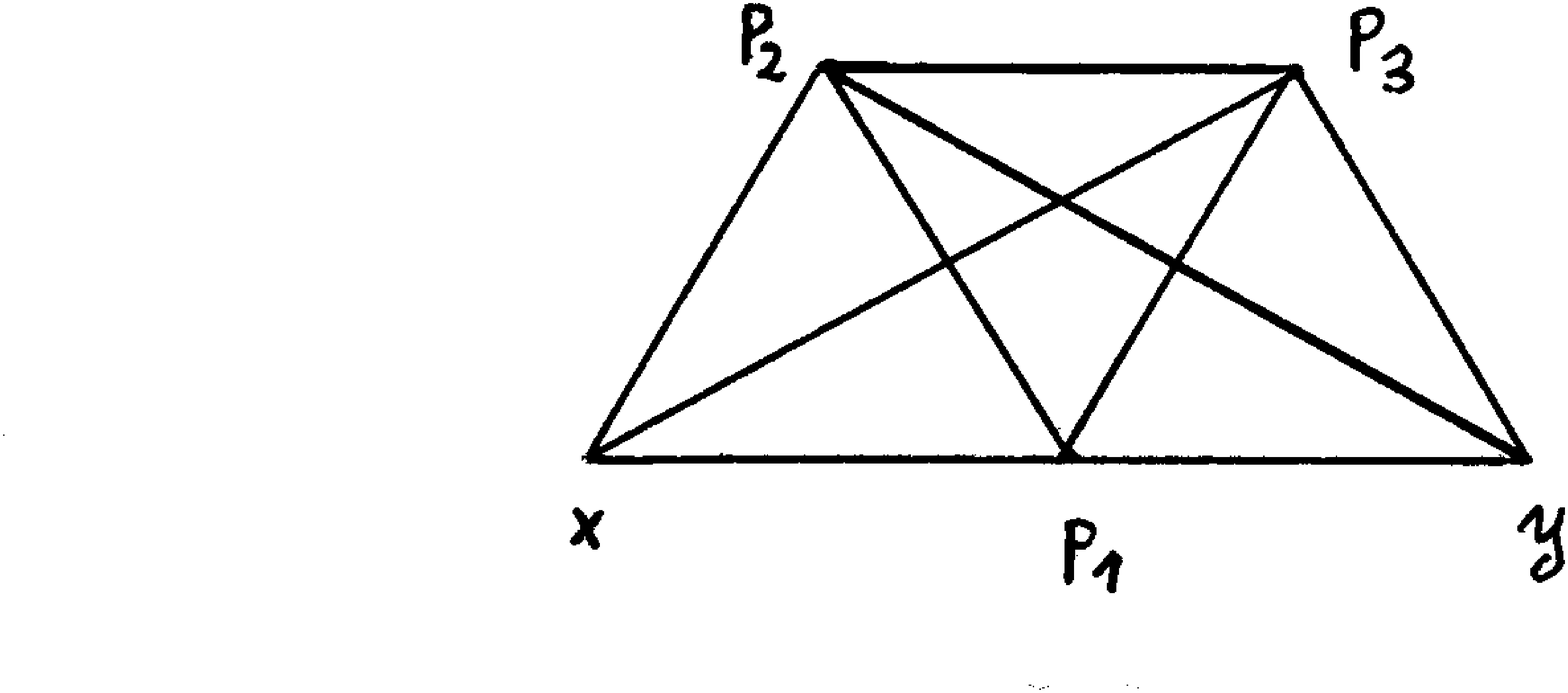}
\\
\centerline{Figure~2}
\centerline{$|x-y|=2 \cdot d$}
\centerline{$|p_1-p_2|=|p_1-p_3|=|p_2-p_3|=|x-p_1|=|x-p_2|=|y-p_1|=|y-p_3|=d$}
\centerline{$|x-p_3|=|y-p_2|=\sqrt{3} \cdot d$}
\vskip 0.2truecm
\par
\noindent
Assume that
$f:S_{xy} \rightarrow {\C}^2$ preserves unit distance. Then 
$f$ preserves all distances between $p_i$ and $p_j$ $(1 \leq i < j \leq 3)$,
$x$ and $p_i$ $(1 \leq i \leq 3)$, $y$ and $p_i$ $(1 \leq i \leq 3)$.
By Proposition~2 the Cayley-Menger determinant of points
$f(x)$, $f(p_{1})$, $f(p_{2})$, $f(p_{3})$, $f(y)$ equals $0$ i.e.
\footnotesize
$$
\det \left[
\begin{array}{ccccccc}
0 &  1                     &  1                       &  1                     &  1                     &  1                        \\
1 & \varphi_2(f(x),f(x))   & \varphi_2(f(x),f(p_1))   & \varphi_2(f(x),f(p_2))   & \varphi_2(f(x),f(p_3))   & \varphi_2(f(x),f(y))  \\
1 & \varphi_2(f(p_1),f(x)) & \varphi_2(f(p_1),f(p_1)) & \varphi_2(f(p_1),f(p_2)) & \varphi_2(f(p_1),f(p_3)) & \varphi_2(f(p_1),f(y))\\
1 & \varphi_2(f(p_2),f(x)) & \varphi_2(f(p_2),f(p_1)) & \varphi_2(f(p_2),f(p_2)) & \varphi_2(f(p_2),f(p_3)) & \varphi_2(f(p_2),f(y))\\
1 & \varphi_2(f(p_3),f(x)) & \varphi_2(f(p_3),f(p_1)) & \varphi_2(f(p_3),f(p_2)) & \varphi_2(f(p_3),f(p_3)) & \varphi_2(f(p_3),f(y))\\
1 & \varphi_2(f(y),f(x))   & \varphi_2(f(y),f(p_1))   & \varphi_2(f(y),f(p_2))   & \varphi_2(f(y),f(p_3))   & \varphi_2(f(y),f(y))  \\
\end{array}
\right]
=0.
$$
\normalsize
Denoting $t=\varphi_2(f(x),f(y))$ we obtain
$$
\det \left[
\begin{array}{ccccccc}
0 &  1   &  1  &  1   &  1   &  1  \\
1 &  0   & d^2 & d^2  & 3d^2 &  t  \\
1 & d^2  &  0  & d^2  & d^2  & d^2 \\
1 & d^2  & d^2 &  0   & d^2  & 3d^2\\
1 & 3d^2 & d^2 & d^2  &  0   & d^2 \\
1 &  t   & d^2 & 3d^2 & d^2  &  0  \\
\end{array}
\right]
=0.
$$
\par
\noindent
Computing this determinant we obtain
$$3d^4 \cdot (t-4d^2)^{2}=0.$$
Therefore
$$t=\varphi_2(f(x),f(y))=\varphi_2(f(y),f(x))=(2d)^2.$$
\vskip 0.2truecm
\par
{\bf Lemma~4.} If $a,b \in D$ and $a>b$, then $\sqrt{a^2-b^2} \in D$.
\vskip 0.2truecm
\par
{\it Proof.}
Let $a, b \in D$, $a>b$, $x,y \in {\R}^2$ and $|x-y|=\sqrt{a^2-b^2}$. Using the notation
of Figure~3 we show that
$$
S_{xy}:=S_{xp_1} \cup S_{xp_2} \cup S_{yp_1} \cup S_{yp_2} \cup S_{p_1p_2}
$$
(where $S_{p_1p_2}$ is known to exist by Lemma~3)
satisfies condition $(\ast)$.
\\
\includegraphics[width=137mm,height=63mm]{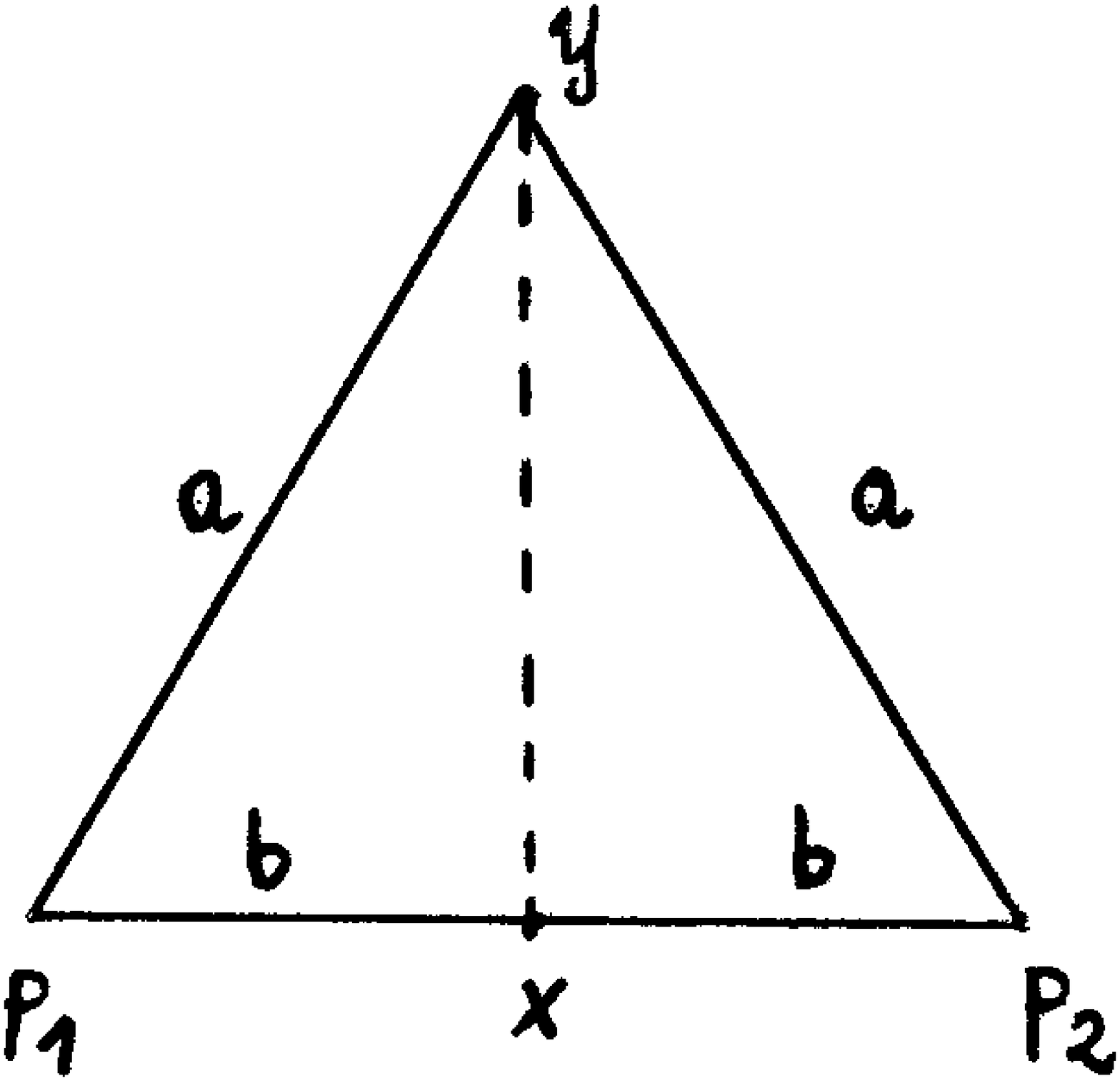}
\\
\centerline{Figure~3}
\centerline{$|x-y|=\sqrt{a^2-b^2}$}
\centerline{$|x-p_1|=|x-p_2|=b$, $|y-p_1|=|y-p_2|=a$, $|p_1-p_2|=2b$}
\vskip 0.2truecm
\par
\noindent
Assume that
$f: S_{xy} \rightarrow {\C}^2$ preserves unit distance.
Then $f$ preserves the distances between
$x$ and $p_{i}$ ($i=1,2$), $y$ and $p_{i}$ ($i=1,2$),
$p_1$ and $p_2$.
By Proposition~2 the Cayley-Menger determinant of points
$f(x)$, $f(p_1)$, $f(p_2)$, $f(y)$ equals $0$ i.e.
\footnotesize
$$
\det \left[
\begin{array}{ccccccc}
0 &  1                     &  1                       &  1                       &  1                    \\
1 & \varphi_2(f(x),f(x))   & \varphi_2(f(x),f(p_1))   & \varphi_2(f(x),f(p_2))   & \varphi_2(f(x),f(y))  \\
1 & \varphi_2(f(p_1),f(x)) & \varphi_2(f(p_1),f(p_1)) & \varphi_2(f(p_1),f(p_2)) & \varphi_2(f(p_1),f(y))\\
1 & \varphi_2(f(p_2),f(x)) & \varphi_2(f(p_2),f(p_1)) & \varphi_2(f(p_2),f(p_2)) & \varphi_2(f(p_2),f(y))\\
1 & \varphi_2(f(y),f(x))   & \varphi_2(f(y),f(p_1))   & \varphi_2(f(y),f(p_2))   & \varphi_2(f(y),f(y))  \\
\end{array}
\right]
=0.
$$
\normalsize
Denoting $t=\varphi_2(f(x),f(y))$ we obtain
$$
\det \left[
\begin{array}{ccccccc}
0 &  1  &  1   &  1   &  1 \\
1 &  0  & b^2  & b^2  &  t \\
1 & b^2 &  0   & 4b^2 & a^2\\
1 & b^2 & 4b^2 &  0   & a^2\\
1 &  t  & a^2  & a^2  &  0 \\
\end{array}
\right]
=0.
$$
\par
\noindent
Computing this determinant we obtain
$$
-8b^2 \cdot (t+b^2-a^2)^2=0.
$$
Therefore $$t=a^2-b^2.$$
\vskip 0.2truecm
\par
{\bf Lemma~5.} If $d \in D$ then $\sqrt{2} \cdot d \in D$.
\par
{\it Proof.} It follows from equality
$\sqrt{2} \cdot d = \sqrt{(\sqrt{3} \cdot d)^2-d^2}$
and Lemmas~1 and~4.
\vskip 0.2truecm
{\bf Lemma~6.} If $d \in D$ then $(2\sqrt{2}/3) \cdot d \in D$.
\vskip 0.2truecm
\par
{\it Proof.}
Let $d\in D$, $x,y \in {\R}^2$, $|x-y|=(2\sqrt{2}/3) \cdot d$.
Using the notation of Figure~4 we show that
$$
S_{xy}:=S_{y\wt{y}} \cup \bigcup_{i=1}^2 S_{xp_i} \cup \bigcup_{i=1}^2 S_{yp_i} \cup S_{p_1p_2}
\cup
\bigcup_{i=1}^2 S_{x\wt{p_i}}
\cup
\bigcup_{i=1}^2 S_{\wt{y}\wt{p_i}}
\cup S_{\wt{p_1}\wt{p_2}}
$$
(where sets corresponding to distances $\sqrt{3} \cdot d$ are known to exist
by Lemma~1, sets corresponding to distances $\sqrt{2} \cdot d$ are known
to exist by Lemma~5, sets corresponding to distances $3 \cdot d$ are known
to exist by Lemma~2) satisfies condition $(\ast)$.
\\
\includegraphics[width=137mm,height=70mm]{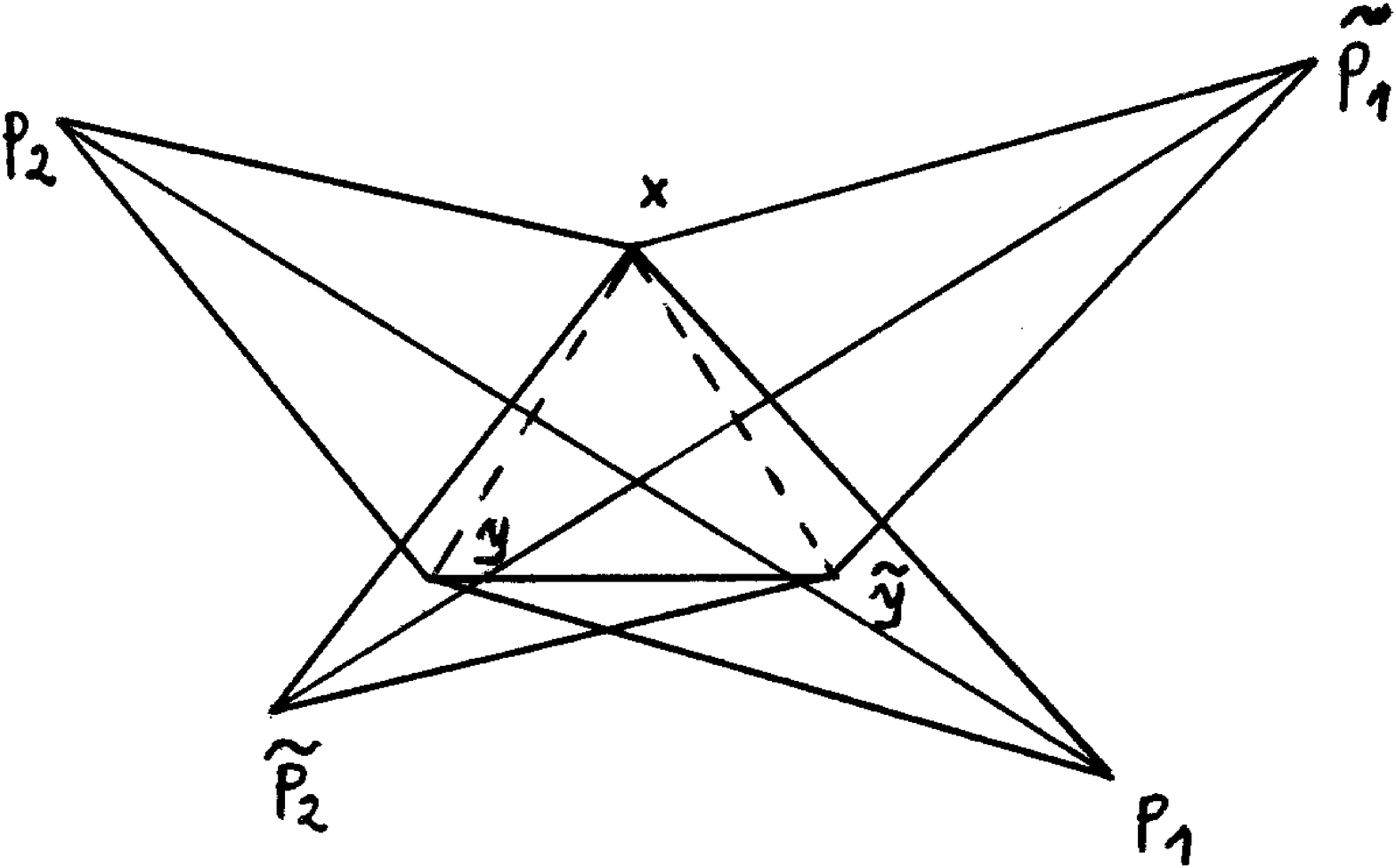}
\\
\centerline{Figure~4}
\centerline{$|x-y|=|x-\wt{y}|=(2\sqrt{2}/3) \cdot d$}
\centerline{$|y-\wt{y}|=d$}
\centerline{$|x-p_1|=|y-p_1|=|x-\wt{p_1}|=|\wt{y}-\wt{p_1}|=\sqrt{3} \cdot d$}
\centerline{$|x-p_2|=|y-p_2|=|x-\wt{p_2}|=|\wt{y}-\wt{p_2}|=\sqrt{2} \cdot d$}
\centerline{$|p_1-p_2|=|\wt{p_1}-\wt{p_2}|=3 \cdot d$}
\vskip 0.2truecm
\par
\noindent
Assume that $f: S_{xy} \rightarrow {\C}^2$
preserves unit distance. Since
$$
S_{xy} \supseteq
S_{y\wt{y}}
\cup
\bigcup_{i=1}^{2}S_{xp_{i}}
\cup
\bigcup_ {i=1}^{2}S_{yp_{i}}
\cup
S_{p_1p_2}
$$
\noindent
we conclude that $f$ preserves the distances between
$y$ and $\wt{y}$,
$x$ and $p_{i}$ $(i=1,2)$, $y$ and $p_i$ $(i=1,2)$,
$p_1$ and $p_2$.
By Proposition~2 the Cayley-Menger determinant of points
$f(x)$, $f(p_1)$, $f(p_2)$, $f(y)$ equals $0$ i.e.
\footnotesize
$$
\det \left[
\begin{array}{ccccccc}
0 &  1                     &  1                       &  1                       &  1                    \\
1 & \varphi_2(f(x),f(x))   & \varphi_2(f(x),f(p_1))   & \varphi_2(f(x),f(p_2))   & \varphi_2(f(x),f(y))  \\
1 & \varphi_2(f(p_1),f(x)) & \varphi_2(f(p_1),f(p_1)) & \varphi_2(f(p_1),f(p_2)) & \varphi_2(f(p_1),f(y) \\
1 & \varphi_2(f(p_2),f(x)) & \varphi_2(f(p_2),f(p_1)) & \varphi_2(f(p_2),f(p_2)) & \varphi_2(f(p_2),f(y))\\
1 & \varphi_2(f(y),f(x))   & \varphi_2(f(y),f(p_1))   & \varphi_2(f(y),f(p_2))   & \varphi_2(f(y),f(y))  \\
\end{array}
\right]
=0.
$$
\normalsize
Denoting $t=\varphi_2(f(x),f(y))$ we obtain
$$
\det \left[
\begin{array}{ccccccc}
0 &  1   &  1   &  1   &  1  \\
1 &  0   & 3d^2 & 2d^2 &  t  \\
1 & 3d^2 &  0   & 9d^2 & 3d^2\\
1 & 2d^2 & 9d^2 &  0   & 2d^2\\
1 &  t   & 3d^2 & 2d^2 &  0  \\
\end{array}
\right]
=0.
$$
\par
\noindent
Computing this determinant we obtain
$$2d^2t \cdot (8d^2-9t)=0.$$
Therefore
$$t=\varphi_2(f(x),f(y))=\varphi_2(f(y),f(x))=((2\sqrt{2}/3) \cdot d)^2$$
or $$t=\varphi_2(f(x),f(y))=\varphi_2(f(y),f(x))=0.$$
Analogously we may prove that
$$\varphi_2(f(x),f(\wt{y}))=\varphi_2(f(\wt{y}),f(x))=
((2\sqrt{2}/3) \cdot d)^2$$
or $$\varphi_2(f(x),f(\wt{y}))=\varphi_2(f(\wt{y}),f(x))=0.$$
If $t=0$ then the points $f(x)$ and $f(y)$ satisfy:
\\
\\
\centerline{$\varphi_2(f(x),f(x))=0=\varphi_2(f(y),f(x))$,}
\\
\centerline{$\varphi_2(f(x),f(p_1))=3d^2=\varphi_2(f(y),f(p_1))$,}
\\
\centerline{$\varphi_2(f(x),f(p_2))=2d^2=\varphi_2(f(y),f(p_2))$.}
\\
\par
\noindent
By Proposition~3b the points
$f(x),f(p_1),f(p_2)$ are affinely independent.
Therefore by Proposition~4 $f(x)=f(y)$ and consequently
$$
d^2=\varphi_2(f(y),f(\wt{y}))=\varphi_2(f(x),f(\wt{y}))
\in \{((2\sqrt{2}/3) \cdot d)^2,~0 \}.
$$
Since $d^2 \neq ((2\sqrt{2}/3) \cdot d)^2$ and $d^2 \neq 0$ we conclude that
the case $t=0$ cannot occur. This completes the proof of Lemma~6.
\vskip 0.2 truecm
\par
Obviously $1 \in D$. Therefore by Lemmas~1 and~6
$$
\{(2\sqrt{2}/3)^k \cdot (\sqrt{3})^l:   k,l \in \{0,1,2,...\}\} \subseteq D.
$$
This completes the proof of Theorem~1.
\vskip 0.2truecm
\par
As a corollary of Theorem~1 and Proposition~5 we obtain our
main theorem.
\vskip 0.2truecm
\par
{\bf Theorem~2.} Each continuous map from ${\R}^2$ to ${\C}^2$ 
preserving unit distance preserves all distances.
\vskip 0.2truecm
\par
{\bf Remark} (\cite{zlozona}). By an endomorphism of $\C$ we understand
any map $f:\C \to \C$ satisfying:
\par
\centerline
{$\forall x,y \in \C ~~f(x+y)=f(x)+f(y)$,}
\centerline{$\forall x,y \in \C ~~f(x \cdot y)=f(x) \cdot f(y)$,}
\centerline{$f(0)=0$,}
\centerline{$f(1)=1$.}
\par
\noindent
If $f:\C \to \C$ is an endomorphism then
$(f_{|\R},...,f_{|\R}): {\R}^n \to {\C}^n$ preserves unit distance.
Bijective endomorphisms are called automorphisms.
There are two trivial automorphisms of $\C$:
identity and conjugation. It is known
that there exist non-trivial automorphisms of $\C$
and each such automorphism $f:\C \to \C$ satisfies:
$\exists^{x \in \R}_{x \neq 0} f(x) \not\in \R$ (\cite{Kuczma}). From this 
$\varphi_n((0,0,...,0),(x,0,...,0))=|x|^2$
and
$\varphi_n((f(0),f(0),...,f(0)),(f(x),f(0),...,f(0)))=(f(x))^2 \neq |x|^2$.
Therefore $(f_{|\R},...,f_{|\R})$ preserves unit distance, but does not
preserve the distance $|x|$~$>$~$0$.

Apoloniusz Tyszka\\
Technical Faculty\\
Hugo Ko{\l}{\l}\c{a}taj University\\
Balicka 104\\
30-149 Krak\'ow, Poland\\
E-mail: {\it rttyszka@cyf-kr.edu.pl}\\
\end{document}